\documentclass[reqno,a4paper]{amsart}

\usepackage{amsfonts,amsmath,amssymb,amsthm}
\usepackage{alltt,esvect,graphicx,stmaryrd,tikz,xcolor}
\usepackage{pgffor}
\usepackage[arrow, matrix, curve]{xy}
\usepackage[titletoc]{appendix}
\usepackage[linesnumbered]{algorithm2e}
\usepackage[T1]{fontenc}
\usepackage[utf8]{inputenc}
\usepackage{helvet}
\usetikzlibrary{calc,intersections,through,shapes,arrows}
\usepackage{mathtools}
\usepackage{subfig}
\usepackage{url}
\usepackage{booktabs}
\usepackage{comment}
\usepackage[textsize=tiny]{todonotes}

\usepackage{enumitem} 
\usetikzlibrary{cd}
\usetikzlibrary{patterns}
\usetikzlibrary{math}
\usetikzlibrary{shapes}
\usetikzlibrary{decorations.fractals}

\usepackage{xstring}

\usepackage{cancel}
\usepackage{hyperref}

\newtheorem{theorem}{Theorem}[section]

\newtheorem{proposition}[theorem]{Proposition}

\theoremstyle{definition}

\newtheorem{example}[theorem]{Example}

\newtheorem{remark}[theorem]{Remark}

\definecolor{colorX}{RGB}{240,0,0}
\definecolor{colorY}{RGB}{0,0,160}
\definecolor{colorZ}{RGB}{0,160,0}

\newcommand{\N}{\mathbb{N}}
\newcommand{\Z}{\mathbb{Z}}
\newcommand{\Q}{\mathbb{Q}}
\newcommand{\R}{\mathbb{R}}

\newcommand{\A}{\mathbb{A}}
\newcommand{\T}{\mathbb{T}}
\newcommand{\V}{\mathbb V}

\renewcommand{\P}{\mathbb{P}}

\newcommand{\CM}{{\mathcal M}}

\newcommand{\CE}{\mathcal{E}}
\newcommand{\CH}{\mathcal{H}}
\newcommand{\CO}{\mathcal{O}}

\newcommand{\Hom}{\mathop{\rm Hom}\nolimits}

\renewcommand{\deg}{\mathop{\rm deg}\nolimits}

\newcommand{\ord}{\operatorname{ord}}

\newcommand{\Spec}{\mr{Spec}\,}

\newcommand{\mr}{\mathrm}

\newcommand{\mc}{\mathcal}

\newcommand{\bp}{\begin{para}}
\newcommand{\ep}{\end{para}}
\newcommand{\bps}{\begin{paras}}
\newcommand{\eps}{\end{paras}}
\newcommand{\benum}{\begin{enumerate}[{\rm(i)}]}
\newcommand{\eenum}{\end{enumerate}}

\newcommand{\Div}{\operatorname{Div}}
\newcommand{\Cl}{\operatorname{Cl}}

\newcommand{\toric}{\T\V}  
\newcommand{\opn}{\operatorname}

\newcommand{\ssect}[1]{Subsection~\ref{#1}}

\newcommand{\Matr}[2]{\left(\begin{array}{@{}*{#1}{r}@{}} #2
\end{array}\right)}

\newcommand{\Matrr}[3]{\left(\begin{array}{@{}*{#1}{r}|*{#2}{r}@{}} #3
\end{array}\right)}

\definecolor{oliwkowy}{HTML}{627037}
\definecolor{lightblue}{RGB}{135,206,250}
\definecolor{darkblue}{RGB}{0,0,160}
\definecolor{darkgreen}{RGB}{0,160,0}
\definecolor{veryPeri}{RGB}{102,103,171}
\definecolor{intOrange}{rgb}{1.0,.310,.0}
\definecolor{MidnightBlue}{RGB}{102,103,171}

\newcounter{para}[section]
\newenvironment*{para}[1]{\refstepcounter{para}\noindent\ignorespaces{\bf\thepara.~#1.}}{\ignorespacesafterend\bigskip}
\newenvironment*{paras}[1]{{\bf #1.}}{\ignorespacesafterend\bigskip}
\numberwithin{para}{section}

\newcommand{\PP}{\mathbb P}

\setlist[enumerate,1]{label = (\roman*),ref = \theenumii.\roman*}

\newcommand{\kk}{{\operatorname{k}}}

\renewcommand{\div}{\operatorname{div}}\newif\ifnotext

\definecolor{colorAB}{RGB}{240,0,0}
\definecolor{colorBC}{RGB}{0,0,160}
\definecolor{colorCD}{rgb}{1.0,.310,.0}
\definecolor{colorDA}{RGB}{0,160,0}

\newcommand{\TT}{{\opn{T}}}

\definecolor{cRayA}{rgb}{0.0,0.6,0.0}
\definecolor{cRayB}{rgb}{0.0,0.0,1.0}
\definecolor{cRayC}{rgb}{1.0,0.0,0.0}

\definecolor{cConeA}{rgb}{0.0,0.6,0.0}
\definecolor{cConeB}{rgb}{0.0,0.0,1.0}
\definecolor{cConeC}{rgb}{1.0,0.0,0.0}

\newif\iftext

\definecolor{green}{rgb}{0.0,0.6,0.0}
\definecolor{blue}{rgb}{0.0,0.0,1.0}
\definecolor{red}{rgb}{1.0,0.0,0.0}

\mathchardef\mhyphen="2D

\newcommand{\OD}{\mathcal O_X(\WD)}
\newcommand{\CV}{\mathcal V}

\newcommand{\idm}{\mathfrak{m}}
\newcommand{\hatDiv}{\operatorname{\Div}}

\newcommand{\ux}{\underline{x}}

\newcommand{\PD}{{P}}

\newcommand{\HM}{{\CH\hspace{-0.09em}\CM}}

\newcommand{\WD}{\mathcal{W}}

\newcommand{\KS}{{K}}
\newcommand{\MN}[2]{{\langle#1,#2\rangle}}

\newcommand{\oldPlusWirdMinus}{{\color{black}-}}

\newif\ifnotext

\begin{document}
\parindent0mm

\title[Horrocks-Mumford sheaves on toric varieties]
{Horrocks-Mumford sheaves on toric varieties}

\author[K.~Altmann]{Klaus Altmann}
\address{Institut f\"ur Mathematik,
FU Berlin,
K\"onigin-Luise-Str.~24-26,
D-14195 Berlin
}
\email{altmann@math.fu-berlin.de}
\author[A.~Hochenegger]{Andreas Hochenegger}
\address{
Dipartimento di Matematica ``Francesco Brioschi'',
Politecnico di Milano,
via Bonardi 9,
20133 Milano 
}
\email{andreas.hochenegger@polimi.it}
\author[F.~Witt]{Frederik Witt}
\address{
Fachbereich Mathematik,
U Stuttgart,
Pfaffenwaldring 57,
D-70569 Stuttgart
}
\email{witt@mathematik.uni-stuttgart.de}

\thanks{{\bf MSC 2020:}
14C20, 
14F06, 
14M25, 
14Q15  
\hfill\newline
\indent{\bf Key words:} computer algebra, Cox ring, divisors, toric geometry, Weil decorations}

\begin{abstract}
Previous work of the authors introduced the class of Horrocks-Mumford sheaves over normal algebraic $\kk$-varieties, the most prominent example being the classical Horrocks-Mumford bundle on projective four-space. For a smooth semiprojective toric variety, we derive an explicit description of these sheaves in terms of finitely generated graded modules over the Cox ring which we use for an algorithmic treatment. In particular, we discuss some explicit examples, their indecomposability and local freeness.
\end{abstract}

\maketitle

\section{Introduction}
\label{sec:Intro}
Let $\kk$ be an algebraically closed field of characteristic $0$ and $X$ a normal $\kk$-variety, that is, a normal, separated and integral scheme of finite type over $\kk$. This enables a well-behaved theory of divisors on $X$ and ultimately of Weil decorations, a notion introduced in~\cite{TRS} and~\cite{WHM}. To define the latter we consider the group of Weil divisors $\Div(X)$ which we endow with the natural poset structure 
\[
D\geq D'\quad:\Longleftrightarrow\quad D-D'\geq0,\text{ that is, }D-D'\text{ is effective.}  
\]
The greatest lower bound or {\em meet} of two divisors is 
\[
D\wedge D':=\min\{D,D'\}:=\sum\min\{D_\PD,D'_\PD\}\cdot\PD, 
\]
where $D_\PD$ is the coefficient of $D\in\Div(X)$ with respect to the prime divisor $\PD$. Given a finite dimensional vector space $\CV$ over the field of rational functions $K=K(X)$, a {\em pre-Weil decoration} is an assignment $\WD\colon\CV\setminus\{0\}\to\hatDiv(X)$ satisfying

\smallskip

\begin{itemize}
\item[{\rm (W1)}] $\WD(f\cdot v)=\div(f)+\WD(v)$ for all $f\in K$ and $v\in\CV$

\smallskip

\item[{\rm (W2)}] $\WD(v+v')\geq\WD(v)\wedge\WD(v')$ for all $v$, $v'\in\CV$;
\end{itemize}
as for ordinary valuations we tacitly extend $\WD$ to all of $\CV$ by sending $0$ to a formal upper bound $\infty$ of $(\Div(X),\leq)$. With $\WD$ we associate the quasi-coherent sheaf
\[
\OD(U):=\{v\in\CV\mid\WD(v)|_U\geq0\}\subseteq\CV
\]
whose generic stalk is $\CV$. We call $\WD$ a {\em Weil decoration}, if $\OD$ is coherent. In this case, $\OD$ is automatically reflexive.

\begin{example}
\label{exam:RefSheaf1}
In view of (W1), a pre-Weil decoration $\WD\colon K\to\hatDiv(X)$ is already determined by $D=\WD(1)$ in $\Div(X)$. Its associated sheaf $\OD$ is precisely the sheaf inside $K$ given by
\[
\OD(U)=\{f\in K^*\mid\big(D+\div(f)\big)|_U\geq0\}\cup\{0\}=\CO_X(D)(U)
\]
on $U\subseteq X$ open. Of course, $\CO_X(\WD)$ is reflexive hence $\WD$ is a Weil decoration.
\end{example}

\begin{remark}
Conversely, any reflexive sheaf $\CE$ gives rise to a Weil decoration $\WD_\CE\colon\CE_\eta\to\Div(X)$ by assigning to $e\in\CE_\eta\setminus\{0\}$ in the generic stalk of $\CE$ the Weil divisor of the reflexive rank one sheaf 
\[
\CE(e)(U):=\big(K\cdot e\big)\cap\CE(U)\stackrel{1/e}{\hookrightarrow}\KS.
\]
In fact, we can boost Weil decorations into a category equivalent to the category of reflexive sheaves, see~\cite
{WHM}.
\end{remark}

For the central example of this paper we consider a smooth and semiprojective toric variety $X=\toric(\Sigma)$ given by the fan $\Sigma$. Further, we let
\begin{equation}
\label{eq:DivT}
\Div_\TT(X):=\bigoplus_{\rho\in\Sigma(1)}\Z D_\rho\subseteq\Div(X) 
\end{equation}
be the subgroup of torus-invariant divsors. This defines a free $\Z$-module with basis provided by the divisors $D_\rho$ corresponding to the one-dimensional cones or {\em rays} $\rho\in\Sigma(1)$; see Section~\ref{sec:Toric} for a brief recap on toric geometry. Let $u\colon\Div_\TT(X)\to\Div_\TT(X)$ be a $\Z$-linear endomorphism which we write as an integer matrix with respect to the basis $\{D_\rho\}_{\rho\in\Sigma(1)}$. We call $u$ {\em admissible} if $u_{\rho\rho}=0$ for all $\rho$ and $[u(D)]=0$ for the class of the divisor $u(D)$; the latter is considered as an element in $M$ with which we can associate the rational function $\chi^{u(D)}$. From~\cite
{WHM} we obtain a Weil decoration $\WD\colon K^2\to\Div(X)$ defined by
\[
\WD_u(f,g)_\PD=
\begin{cases}
\min\{\ord_P(f),\ord_P(g)\}+1\!\!&\text{if }\PD=D_\rho\text{ and }
(f/g)(D_\rho)=\chi^{u(D_\rho)}\\
\min\{\ord_P(f),\ord_P(g)\}\!\!&\text{else}.
\end{cases}
\]
Here, $\ord_\PD$ denotes the valuation of the local ring $\CO_{X,\PD}$, $\PD$ a prime divisor, and evaluation at $D_\rho$ yields an element $(f/g)(D_\rho)$ in the residue field $\CO_{X,D_\rho}/\idm_{X,D_\rho}$; since the diagonal of $u$ vanishes, we may consider $\chi^{u(D_\rho)}$ as a nontrivial element of this field.

\begin{example}
\label{exam:HMType}
If we let $D_0,\ldots,D_4$ be the divisors corresponding to the rays $\rho_0,\ldots,\rho_4$ of $\PP^4$ (see Example~\ref{exam:PnToric}), then 
in~\cite
{WHM} we show that the admissible integer matrix
\begin{equation}
\label{eq:MatHM}
u=\Matr{5}{ 
0&1&-1&-1&1\\
1&0&1&-1&-1\\
-1&1&0&1&-1\\
-1&-1&1&0&1\\
1&-1&-1&1&0}
\end{equation}
defines the celebrated Horrocks-Mumford bundle~\cite{hm}.
\end{example}

We refer to the reflexive sheaf associated with $\WD_u$ as {\em H(orrocks)-M(umford) sheaf} and denote it $\HM(X,u)$ (they can be defined in greater generality, see~\cite{WHM}). Our main theoretical result for these sheaves is this.

\medskip

{\bf Theorem (see~\ref{thm:GradMod}).} 
{\em Let $X$ be a smooth and semiprojective toric variety with Cox ring $S(X)$, and $u$ be an admissible $r\times r$-integer matrix with $r=\sharp\Sigma(1)$. Then we construct an explicit finitely generated and graded $S(X)$-module $G(X,u)$ whose sheafification yields $\HM(X,u)$.}

\medskip

This explicit representation makes $\HM$-sheaves over toric varieties amenable to algorithmic methods with $G(X,u)$ as input. We illustrate this by a series of concrete examples in Section~\ref{sec:AlgAsp} for which we compute Chern classes and discuss irreducibility and local freeness, based on our {\tt Julia}-package~\cite{WD} for the computer algebra system~\cite{OSCAR}. 

\bigskip

{\em Acknowledgements.} We thank Lars Kastner and Max Horn for their valuable support with {\tt Oscar}.

\bigskip

{\em Declaration of AI-assisted technologies.}
During the preparation of this article, the authors used AI tools for coding in Oscar. All AI-assisted code was reviewed, edited, and independently checked by the authors.

\section{Toric geometry}
\label{sec:Toric}
From now on, we let $X=\toric(\Sigma)$ be an $n$-dimensional toric $\kk$-variety defined by the fan $\Sigma$ and with embedded open torus $\TT$. The {\em character lattice} is the free $\Z$-module given by algebraic group morphisms $\TT\to\kk^*$, that is,
\[
M=\Hom_{\opn{ag}}(\TT,\kk^*)\cong\Z^n.
\]
It induces the $\kk$-algebra $\kk[M]$ for which $\Spec\kk[M]=\TT$. More generally, every cone $\sigma\in\Sigma$ gives rise to a toric variety $\toric(\sigma)$ which we can equivariantly embed as the open affine subset 
\[
U_\sigma=\Spec\kk[\sigma^\vee\cap M]. 
\]
Toric geometry is thus based on the interplay between geometry, algebra and the combinatorics of the fan. For instance, $X$ is smooth if and only if the minimal generators of any cone form part of a $\Z$-basis of $N$. A general reference for toric varieties is~\cite{CoxBook}.

\medskip

For technical reasons, we discard some degenerate cases and {\em always} assume that the finite set of one-dimensional cones or {\em rays} $\Sigma(1)$ generates $N_\R=N\otimes_\Z\R$, where
\[
N=\Hom_\Z(M,\Z) 
\]
is the dual of $M$. Then the fundamental sequence of toric geometry reads as
\begin{equation}
\label{eq:FundSeq}
\xymatrix{
0\ar[r]&M\ar[r]^-{\div}&\Div_\TT(X)=\bigoplus\limits_{\rho\in\Sigma(1)}^{\phantom{a}}\Z D_\rho\ar[r]^-{[\cdot]}&\Cl(X)\ar[r]&0,
}
\end{equation}
cf.~\eqref{eq:DivT}. Here, an element $m$ in $M$ gets mapped to the divisor of the regular function $\chi^m$ on $\TT$ defined by $m\in M$. If $M\times N\to\Z$, $(m,n)\mapsto\MN{m}{n}$ denotes the natural pairing, then
\[
\div(\chi^m)=\sum_{\rho\in\Sigma(1)}\MN{m}{\rho}D_\rho. 
\]
Second, $[\cdot]$ sends a toric divisor to its class. As usual, we shall identify a ray $\rho$ in $\Sigma(1)$ with its primitive generator.

\begin{example}
\label{exam:PnToric}
For $X=\PP^n$ the fan is generated by the rays
\[
\rho_0=-\sum_{i=1}^ne_i,\;\rho_1=e_1,\;\ldots,\;\rho_n=e_n 
\]
for the standard basis $e_1,\ldots,e_n$ of $M=\Z^n$.
\end{example}

In general, we may think of toric varieties as a generalisation of projective spaces via the so-called {\em Cox construction}. Given $X=\toric(\Sigma)$, its {\em Cox ring} is
\[
S:=\kk[x_\rho\mid\rho\in\Sigma(1)]. 
\]
The {\em degree} of an $S$-monomial $x^a=\prod_{\rho\in\Sigma(1)}x_\rho^{a_\rho}$, $a=(a_\rho)\in\Z^{\Sigma(1)}$, is defined to be
\[
\deg x^a:=\big[\sum a_\rho D_\rho\big]\in\opn{Cl}(X). 
\]
This turns $S$ into a graded $\kk$-algebra and gives rise to the local coordinate rings by homogeneous localisation. Indeed, the localisation at \fbox{$x_{\hat\sigma}=\prod_{\rho\not\in\sigma(1)}x_\rho$} is given by
\[
S_{x_{\hat\sigma}}=\kk[x_\rho,\,x^{-1}_\tau\mid\rho\in\Sigma(1),\,\tau\not\in\sigma(1)]
\]
and inherits the natural $\opn{Cl}(X)$-grading above. In particular, the degree $0$-part is
\begin{equation}
\label{eq:HomLoc}
S_{(x_{\hat\sigma})}:=(S_{x_{\hat\sigma}})_0=\kk[\sigma^\vee\cap M]
\end{equation}
where the last equality follows from identifying $\chi^m\in\kk[\sigma^\vee\cap M]$ with
\[
\underline x^m:=\prod_{\rho\in\Sigma(1)}x_\rho^{\langle m,\rho\rangle};
\]
note that $\div(\chi^m)$ is indeed of $\opn{Cl}(X)$-degree $0$ by the Fundamental Sequence~\eqref{eq:FundSeq}. More generally, any coherent sheaf $\mc F$ can be obtained in this way, viz., there exists a (not necessarily unique) finitely generated graded $S$-module $F$ such that
\[
\Gamma\big(U_\sigma,\mc F\big)=F_{(x_{\hat\sigma})} 
\]
for any $\sigma\in\Sigma(n)$.

\section{$\HM$-sheaves over toric varieties}
\label{sec:HMSheaves}
We now turn to the sheaves $\HM(X,u)$ from Section~\ref{sec:Intro}\,, where $X=\toric(\Sigma)$ is a smooth and semiprojective toric variety, viz., $\Sigma$ is the normal fan of a full-dimensional (not necessarily bounded) lattice polyhedron in $M$, and $u=(u_{\rho\rho'})$ is an admissible integer $r\times r$-matrix with $r=\sharp\Sigma(1)$. As a result, we can cover $X$ by the open affines
\[
U_\sigma=\Spec\kk[\sigma^\vee\cap M]=\Spec\kk[x_\rho\!\mid\!\rho\in\sigma(1)]=\A^n,\quad\sigma\in\Sigma(n).
\]
If we let $u_{\rho'}$, $\rho'\in\Sigma(1)$, be the element in $M$ given by the $\rho'$-th column of $u$, then
\[
\kk[M]\ni\chi^{u_\rho'}=\prod\limits_{\rho\in\Sigma(1)}x^{\langle u_{\rho'},\rho\rangle}_\rho=\prod\limits_{\rho\in\Sigma(1)}x^{u_{\rho\rho'}}_\rho=\underline x^{u_\rho'}.
\]

\subsection{The trivial $\HM$-sheaves \texorpdfstring{$\HM$}{HM}$(X,0)$}
For $X=\toric(\Sigma)$ we denote
\[
\partial X:=\sum_{\rho\in\Sigma(1)}D_\rho  
\]
its boundary divisor.

\begin{proposition}
We have an explicit isomorphism $\HM(X,0)=\CO_X\oplus\CO_X(\partial X)$.
\end{proposition}

\begin{proof}
Consider the complex 
\[
\xymatrix{
0\ar[r]&\HM(X,0)\ar[r]&\CO_X(\partial X)^2\ar[r]^-\Psi&\CO_X\big(\partial X\big)\Big/\CO_X\ar[r]&0 
}
\]
induced by $\Psi(f,g)=\overline{f-g}$, the residue class of $f\oldPlusWirdMinus g$. Since $X$ is smooth and semiprojective, the torus-invariant open affines $U_\sigma=\Spec\kk[x_\rho\,|\,\rho\in\sigma(1)]=\A^n$ associated with the top-dimensional cones $\sigma\in\Sigma(n)$ cover $X$. The local bases $e_0=(1,1)/\prod_{\rho\in\sigma(1)}\!x_\rho$ and $e_1=(1,0)$ of the kernel of $\Psi$ glue to the diagonal $\Delta\big(\CO_X(\partial X)\big)$ and $\CO_X\oplus0$ inside $\CO_X(\partial X)^2$.
\end{proof}

\subsection{\texorpdfstring{$\HM$}{HM}-sheaves on $X=\A^n$}
\label{subsec:HMAff}
Fix coordinates $x=(x_1,\ldots,x_n)$ on $\A^n$ so that
\[
M=\Z^n\quad\text{and}\quad\Div_\TT(\A^n)=\bigoplus_{i=1}^n\Z H_i
\]
for the coordinate hyperplanes $H_i=\{x_i=0\}$. The sheaf $\HM(u):=\HM(\A^n,u)$ is specified by an $n\times n$ integer matrix $u$ with trivial diagonal. Let \fbox{$G(u)=G(\A^n,u)$} be the corresponding reflexive $\kk[\ux]$-module of global sections $\Gamma\big(\A^n,\HM(u)\big)$. We say that
\begin{equation}
\label{eq:nuCond}
(f,g)\in K^2\text{ satisfies }(\ast)_i\quad\text{if and only if}\quad 
f/g\equiv\ux^{u(H_i)}\opn{mod}x_i\kk[\ux]_{(x_i)}\,.
\end{equation}
Further, $p:=x_1\cdot\ldots\cdot x_n$ and $p_i:=p/x_i$ for $i=1,\ldots,n$, and for $U_i:=[x_i\neq 0]=\A^n\setminus H_i$ we set
\[
V_i:=\toric(\rho_i)=\{p_i\neq0\}=\bigcap_{j\neq i} U_j, 
\]
that is, $V_i$ is the union of the torus $\TT^n$ of $\A^n$ as well as of the inner points of the toric divisor $H_i$. The localisation $\kk[\ux]_{p_i}$ represent the regular functions on $V_i$ and
\[
G(u)_{p_i}=\Gamma(V_i,\HM(u))=\{(f,g)\in\frac{1}{x_i}\cdot\kk[\ux]_{p_i}^2\mid(f,g)\text{ satisfies }(\ast)_i\text{ or }f,\,g\in\kk[x]_{p_i}\}.
\]
If $u^+(i)$ and $u^-(i)\in\N^n$ denote the positive and negative part of the column $u(H_i)\in\Z^n$, that is, $u(H_i)=u^+(i)-u^-(i)$, condition $(\ast)_i$ reads precisely as
\[
f\in\frac{1}{x_i}\cdot\underline x^{u^+(i)}\cdot h(x)+\big(\kk[\ux]_{p_i}\big)_{(x_i)}\quad\text{and}\quad g\in\frac{1}{x_i}\cdot\underline x^{u^-(i)}\cdot h(x)+\big(\kk[\ux]_{p_i}\big)_{(x_i)}
\]
for a common function $h(x)\in\kk[\ux]_{p_i}$ which doesn't contain the 
variable $x_i$. Therefore, the $\kk[\ux]_{p_i}$-module $G(u)_{p_i}$ 
is generated inside $\frac{1}{x_i}\cdot\kk[\ux]_{p_i}^2$ by the basis
\begin{equation}
\label{eq:BasisVi}
f_1=(x_i,0)/x_i\quad\text{and}\quad f_2=(\underline x^{u^+(i)},\,\underline x^{u^-(i)})/x_i.
\end{equation}
The $\kk[x]$-module $G(u)$ on $\A^n$ can now be computed via the following two methods.

\subsubsection{Exhausting \texorpdfstring{$G(u)$}{G(u)}:~Taking the reflexive hull}
\label{subsubsec:RefHull}
Consider the $\kk[\ux]$-module
\[
G'(u):=\Big\langle\frac{1}{x_i}\cdot\big(\underline x^{u^+(i)},\, \underline x^{u^-(i)}\big)\mid i=1,\ldots,n\Big\rangle+\kk[\ux]^2\;\subseteq\;\frac{1}{p}\cdot\kk[\ux]^2.
\]
The inclusion $G'(u)\subseteq G(u)$ holds by design. Further,~\eqref{eq:BasisVi} entails $G'(u)_{p_i}=G(u)_{p_i}$ for all $i=1,\ldots,n$. As a result,
\[
G(u)\,\subseteq\,G(u)_{p_i}=G'(u)_{p_i}
\]
whence $G'(u)\subseteq\,G(u)\,\subseteq\,\Gamma(\bigcup_i\!V_i,\HM'(u))$ where $\HM'(u)$ is the sheafification of $G'(u)$. Now the open set 
\[
\bigcup_{i=1}^nV_i\;=\;\toric\big(\Sigma(1)\big)\;=\;\bigcap_{i<j}(U_i\cup U_j) \;=\;
\A^n\setminus\,\bigcup_{i<j}(H_i\cap H_j)
\]
is obtained by removing codimension $2$-subsets of $\A^n$. Reflexivitiy of $G(u)$ finally implies
\[
G(u)\,=\,G(u)^{\vee\vee}\,=\,G'(u)^{\vee\vee},
\]
which is the reflexive hull of the explicitly known module $G'(u)$.

\subsubsection{Envelopping \texorpdfstring{$\HM(u)$}{HM(u)}:~Taking the intersection}
\label{subsubsec:TakeInt}
Turning the tables we consider the bigger $\kk[\ux]$-modules
\begin{equation}
\label{eq:HMI}
G_i(u):=\big\{(f,g)\in\frac{1}{p}\cdot\kk[\ux]^2\mid(f,g)\text{ satisfies }(\ast)_i\big\}
\end{equation}
for $i=1,\ldots,n$. Now pairs of rational functions in $G_i(u)$ still satisfy condition $(\ast)_i$, but may have poles of order one along the divisors $H_j$, $j\neq i$, without any further constraints. Argueing as for~\eqref{eq:BasisVi} shows that $G_i(u)$ is generated by
\begin{equation}
\label{eq:Giu}
G_i(u)\;=\;\big\langle\frac{1}{p}\cdot(\underline x^{u^+(i)},\,\underline x^{u^-(i)}),\,\frac{1}{p}\cdot(x_i,0),\,\frac{1}{p}\cdot(0,x_i)\big\rangle
\end{equation}
so that intersecting in $\frac{1}{p}\cdot\kk[\ux]^2$ yields
\[
G(u)\;=\;\bigcap\limits_{i=1}^nG_i(u)\;=\;\frac 1p\bigcap\limits_{i=1}^n\big\langle(\underline x^{u^+(i)},\,\underline x^{u^-(i)}),\;(x_i,0),\;(0,x_i)\big\rangle 
\]
by design of $G_i(u)$.

\subsection{The first Chern class}
Building upon the previous discussion we will derive two global statements for general $\HM(X,u)$. For the first statement, we recall that two vectors $v_1$, $v_2$ in $K^2$ are $\PD$-{\em orthogonal} if
\[
\WD_u(f_1v_1+f_2v_2)_\PD=\min\{\ord_\PD(f_1),\,\ord_\PD(f_2)\}.
\]
This holds precisely if $(v_1,v_2)$ defines a basis of $\HM(X,u)_\PD$, cf.~\cite
{WHM}.

\begin{proposition}
\label{prop:Chern1}
Let $X$ be a smooth and semiprojective toric variety. Then
\[
\det\HM(X,u)=\CO_X(\partial X).
\]
In particular, we can identify the first Chern class $c_1(\HM(X,u))$ with the image of the anti-canonical divisor in the Chow group $A^1$.
\end{proposition}

\begin{proof}
By definition, the line bundle $\det\HM(X,u)=\big(\Lambda^2\HM(X,u)\big)^{\vee\vee}$ is the reflexive hull of the second exterior power $\Lambda^2\HM(X,u)$. From Example~\ref{exam:RefSheaf1} its Weil decoration is determined by any $0\neq e\in\det\HM(X,u)_\eta$, for instance by $e_1\wedge e_2$ for the standard basis $(e_1,e_2)$ of $K^2=\HM(X,u)_\eta$. This is $\PD$-orthogonal for all prime divisors in $\TT$ whence
\[
\WD_{\det\HM(X,u)}(e_1\wedge e_2)_\PD=0,\quad\PD\in\TT.
\]
For the boundary divisors $D_\rho=\{x_\rho=0\}$ in $X\setminus\TT$, $\rho\in\Sigma(1)$, we restrict $\HM(X,u)$ to $V_\rho=\toric(\rho)$ inside some affine chart given by $\sigma\in\Sigma(n)$ with $\rho\in\sigma$. Further, we consider from~\eqref{eq:BasisVi} the $D_\rho$-orthogonal basis $(f_1,f_2)$ of the stalk $\HM(X,u)_{D_\rho}$ at (the generic point of) $D_\rho$. Then $f_1\wedge f_2$ is $D_\rho$-orthogonal and $e_1\wedge e_2=v_\rho\cdot x_\rho\cdot f_1\wedge f_2$ for some unit $v_\rho$ in $\CO_{X,D_\rho}$. As a result, 
\[
\WD_{\det\HM(X,u)}(e_1\wedge e_2)_{D_\rho}=1,\quad\rho\in\Sigma(1)
\]
whence $\WD_{\det\HM(X,u)}(e_1\wedge e_2)=\partial X$.
\end{proof}

\subsection{\texorpdfstring{$\HM$}{HM}-sheaves as graded modules over the Cox ring}
\label{subsec:HMPPn}
As a coherent sheaf, $\HM(X,u)$ can be obtained from a finitely generated graded $S$-module $G(X,u)$ by sheafification. Using~\ssect{subsec:HMAff}\, for charts as well as for ``Cox cones'' -- the total coordinate spaces -- such a module can be computed as follows.

\begin{theorem}
\label{thm:GradMod}
The sheafification of $G(X,u)=G(\A^{\Sigma(1)}\!,u)$ yields $\HM(X,u)$. Furthermore, for all $\sigma\in\Sigma(n)$ we have
\[
G(U_\sigma,u)=G(\A^n,u_\sigma)=G(\A^{\Sigma(1)},u)_{(x_{\hat\sigma})},
\]
where $u_\sigma$ is the matrix obtained by deleting the $\rho$-th column and row for any $\rho\not\in\sigma$.
\end{theorem}

\begin{proof}
For a top-dimensional cone $\sigma\in\Sigma(n)$ we have the toric divisors $\{D_\rho\}_{\rho\in\sigma(1)}$ on $U_\sigma$, and the datum $u_\rho$ induces 
the linear map $u_\sigma\colon M\to\Z^{\sigma(1)}$ with matrix $u_\sigma=(\langle u_{\rho'},\rho\rangle)_{\rho,\rho'\in\sigma(1)}$ as asserted. The resulting module $G(\A^n,u_\sigma)=\bigcap_{\rho\in\sigma(1)}G_\rho(u)$ with $G_\rho(u)$ as in~\eqref{eq:Giu} is just the localisation of $G(\A^{\Sigma(1)}\!,u)_{(x_{\hat\sigma})}$, where we consider $u$ as defining an $\HM$-sheaf on the total coordinate space $\A^{\Sigma(1)}$ and whose module $G(\A^{\Sigma(1)}\!,u)$ is given by the intersection $\bigcap\limits_{\rho\in\Sigma(1)}G_\rho(u)$.
\end{proof}

\section{Algorithmic aspects}
\label{sec:AlgAsp}
The explicit construction in Theorem~\ref{thm:GradMod} of the graded module $G(X,u)$ enables a computational approach. In the following we make use of the {\tt Julia} based computer algebra system {\tt Oscar}~\cite{OSCAR}; the latter provides, among other things, a symbolic treatment of toric varieties.

\subsection{The {\tt HorrocksMumford} package \cite{WD}}
Let $X$ be a smooth and semiprojective toric variety $X$ and $u$ an admissible matrix. Then the package computes for the sheaf $\HM(X,u)$
\begin{itemize}
    \item generators of the graded $S$-module $G=G(X,u)$ from Theorem~\ref{thm:GradMod} via
    \[
    {\tt G = HM(X,u).}
    \]
    \item its Chern classes
    \[
    {\tt c = chern(X,u)} 
    \]
    provided $X$ is projective. This uses {\tt Oscar}'s command {\tt free\_resolution} for the module $G$. The direct summands of the resolution correspond to invertible sheaves whose Chern classes compute those of $\HM(X,u)$ in the Chow ring of $X$, the latter being again provided by {\tt Oscar} via {\tt chow\_ring(X)}.
\end{itemize}

\begin{remark}
The module $G(X,u)=G(\A^{\Sigma(1)},u)$ can be computed either by the exhausting method from~\ssect{subsubsec:RefHull} or by taking intersections as in \ssect{subsubsec:TakeInt}\,. The method can be imposed via {\tt HM(X,u;method=:refhull)} or {\tt HM(X,u;method=:intersec)}, respectively; default is computation by intersection.
\end{remark}

Further, we can check whether or not
\begin{itemize}
    \item for projective $X$, the total Chern class $\HM(X,u)$ factors in the Chow ring via
    \[
    {\tt is\_num\_indecomposable\_HM(X,u)} 
    \]
    If this is true, that is, the total Chern class is irreducible, then $\HM(X,u)$ must be indecomposable. 
    \item $\HM(X,u)$ is locally free via
    \[
    {\tt is\_locally\_free\_HM(X,u)}; 
    \]
    this uses {\tt Oscar}'s built in function for computing Fitting ideals.
\end{itemize}
See~\cite{WD} for further documentation.

\subsection{Examples}

\begin{example}
\label{exam:P2}
Consider the matrix
\[
u= 
\left(\begin{array}{rrr}
0&-1&1\\
1&0&-1\\
-1&1&0
\end{array}\right)
\]
whose rows add to zero so that $u$ defines the Horrocks-Mumford sheaf $\HM(\P^2,u)$. Over the Cox ring $S=\kk[x_1,x_2,x_3]$, {\tt HorrocksMumford} computes the generators
\[
(-x_2^2x_3,-x_2x_3^3),\ (x_1x_2x_3,0),\ (-x_1^2x_2,-x_1x_2^2),\ (0,-x_1x_2x_3),\ (-x_1x_3^2,-x_1^2x_3)
\]
of the $S$-module $G(\PP^2,u)$ representing $\HM(\PP^2,u)$. The Chern polynomial computed by {\tt Oscar} can be written as
\[
c(\HM(\PP^2,u))=6h^2+3h+1
\]
in the Chow ring $A^\bullet(\P^2) = \Z[h]/h^3$; in particular, $c_1=c_1(\CO_{\PP^2}(3))$ in accordance with Proposition~\ref{prop:Chern1}. Since $1+3h+6h^2$ cannot be written as a product of the form $(1+ah)\cdot(1+bh)$, the rank two bundle $\HM_{\P^2}(u)$ must be indecomposable, which is confirmed by {\tt is\_num\_indecomposable\_HM}. Finally, $\HM_{\P^2}(u)$ is locally free as we can check by {\tt is\_locally\_free\_HM}.
\end{example}

\begin{remark}
Being reflexive, $\HM$-sheaves are locally free on smooth surfaces. First counterexamples exist already in dimension $3$~\cite{WHM}: $\HM(\A^3,u)$ with matrix 
\[
u=\Matr{3}{
0&1&0\\
1&0&-1\\
0&-1&0}
\]
is not free. 
\end{remark}

\begin{example}
The classical Horrocks-Mumford bundle $\HM(\PP^4,u)$ with $u$ given by~\eqref{eq:MatHM} is locally free and indecomposable as can be checked with {\tt HorocksMumford}. In addition, we obtain 19 generators for the corresponding $S$-module $G(\PP^4,u)$. The Chern polynomial is computed as
\[
c(\PP^4,u)=10h^2+5h+1
\]
in accordance with the literature~\cite{hm}
\end{example} 

\begin{example}
As a last example, we consider the smooth toric varieties of Picard rank $2$ given by the total space of the projectivised bundle
\[
X_{n,\underline a}=\PP(\CO_{\PP^n}\oplus\bigoplus_{i=1}^r\CO_{\PP^n}(a_i))\to\PP^n
\]
for an integer vector $\underline a=(a_1,\ldots,a_r)$ with $0\leq a_1\leq\ldots\leq a_r$. There is a natural ordering for the $n+1+r+1$ rays such that the rays $e_1,\ldots,e_{n+1}$ represent the base, while $f_1,\ldots,f_{r+1}$ represent the fibre; the map $\Div_\TT(X)=\Z^{n+1+r+1}\to\opn{Cl}(X)$ is then given by the matrix
\[
\Matrr{3}{4}{1&\ldots&1&0&a_1&\ldots&a_r\\
0&\ldots&0&1&1&\ldots&1}\colon\Z^{n+1+r+1}\to\Z^2,
\]
see~\cite{kle88}. In particular, any admissible matrix $u$ of $\PP^n$ lifts to an admissible matrix $\tilde u$ of $X_{n,\underline a}$ by extending $u$ trivially. We can implement these varieties via the function {\tt kleinschmidt} of the {\tt HorrocksMumford}-package; {\tt Oscar}'s built-in command 
\begin{verbatim}
map_from_torusinvariant_weil_divisor_group_to_class_group
\end{verbatim}
yields the matrix with respect to the rays given by {\tt rays(X)}; note that {\tt Oscar}'s ordering will in general be different from ours. 

\medskip

For instance, consider $X=X_{2,(2)}$ and the matrix $u$ from Example~\ref{exam:P2}. This yields the indecomposable, but not locally free sheaf $\HM(X,\tilde u)$ for
\[
\tilde u=\Matr{5}{
0&-1&1&0&0\\
1&0&-1&0&0\\
0&0&0&0&0\\
0&0&0&0&0\\
-1&1&0&0&0}.
\]
On the other hand,
\[
v=\Matr{5}{
0&-1&-2&1&1\\
-1&0&1&0&-3\\
-1&1&0&0&1\\
1&-1&0&0&-1\\
3&-1&1&-1&0}
\]
yields the indecomposable and locally free sheaf $\HM(X,v)$. For the Chern classes, we first compute the Chow ring via {\tt Oscar} which yields
\[
A=\Q[x_1,\ldots,x_5]/(x_1-x_5,x_2-x_5,-x_3+x_4+2x_5,x_3x_4,x_1x_2x_5)). 
\]
The classes of $x_1$, $x_2$ and $x_5$ are represented by the pullback of a hyperplane section class $h$ of $\PP^2$. Relabeling $x_4=H$ where $H$ is the relative hyperplane section class of $X$, we have $A=\Q[h,H]/(h(h+2H),H^3)$. For $\HM(X,\tilde u)$, {\tt chern} computes the Chern classes
\[
c_1=5h+2H,\quad c_2=6hH+13h^2=10h^2\quad\text{and}\quad c_3=12h^2H 
\]
in accordance with Proposition~\ref{prop:Chern1}; note that $\HM(X,\tilde u)$ is merely reflexive and not locally free so that $c_3$ can indeed be nontrivial. On the other hand, we find
\[
c_1=5h+2H,\quad c_2=12h^2+7hH=h(9h+H)\quad\text{and}\quad c_3=0
\]
for $\HM(X,v)$.
\end{example}

\end{document}